\RequirePackage{fix-cm}
\documentclass[smallextended]{svjour3}       
\smartqed  

\usepackage{graphicx}
\usepackage{amssymb}
\usepackage{amsthm}
\usepackage[fleqn]{amsmath}
\usepackage{mathtools}
\usepackage{cases}
\usepackage{cite}
\usepackage{lipsum}
\usepackage{mathptmx}

\usepackage{latexsym}
\journalname{Nonlinear Dynamics}
\begin{document}

\title{A critical analysis of the conformable derivative}


\author{Ahmed A. Abdelhakim \footnote{Corresponding author}\and
         Jos\'{e} A. Tenreiro Machado}

\authorrunning{Abdelhakim and Machado} 

\institute{Ahmed A. Abdelhakim \at
             Mathematics Department, Faculty of Science, Assiut University, 71516 Assiut - Egypt \\
              \email{ahmed.abdelhakim@aun.edu.eg}           
           \and
           Jos\'{e} A. Tenreiro Machado\at
             Institute of Engineering of Polytechnic of Porto, Department of Electrical Engineering, Porto - Portugal \\
              \email{jtm@isep.ipp.pt}}

\date{Received: date / Accepted: date}
\maketitle
\begin{abstract}
We prove that conformable ``fractional"  differentiability
of a function $f:[0,\infty[\,\longrightarrow \mathbb{R}$ is nothing else than the classical differentiability. More precisely,
the conformable $\alpha$-derivative of $f$ at some point $x>0$, where $0<\alpha<1$, is the pointwise product
$x^{1-\alpha}f^{\prime}(x)$.
This proves the lack of significance of recent studies of the conformable derivatives.
The results imply that
interpreting fractional derivatives
in the conformable sense
alters fractional differential problems
into differential problems with the usual integer-order derivatives that no longer describe the original fractional physical phenomena.
A general fractional viscoelasticity model
is analysed to illustrate this state of affairs.
We also test the modelling efficiency of the conformable derivative using a fractional model of
viscoelastic deformation of tight sandstone, and a fractional world population growth model.
\end{abstract}

\section{Introduction and remarks}
Despite the large volume of research on fractional calculus, no specific criteria for defining a fractional derivative are usually formulated.
Various definitions for fractional differential operators have been reviewed (see e.g., \cite{Giusti1,Kilbas,Oliveira,
Ortigueira,Ortigueira1,Ortigueira2,
Tarasov,Tarasov1,Tarasov2}). The concept of the conformable derivative, on the other hand, has not been scrutinized enough.\\
\indent
It is proved in \cite{Ortigueira} that the Gr\"{u}nwald-Letnikov, Riemann-Liouville and Caputo fractional differential operators share a set of properties that may be expected from a differential operator to be considered fractional. These properties are:
\begin{description}
  \item[(i)] Linearity
  \item[(ii)] The zero order operator is the identity
  \item[(iii)] Integer order operators give the ordinary derivative
  \item[(iv)] Index law and commutativity
  \item[(v)] The generalized Leibniz rule for the derivative of a product:
      \begin{equation*}
D^{\alpha}\left(f(x)g(x)\right)
=\sum_{n=0}^{\infty}\binom{\alpha}{n}
D^{n}f(x)\,D^{\alpha-n}g(x).
      \end{equation*}
\end{description}
\indent In an attempt to define a fractional derivative
that inherits as many features as possible from the classical derivative, Khalil et al. \cite{khalil} proposed defining the ``conformable fractional derivative" $T_{\alpha}f$, for $0<\alpha<1$,  of a function $f:[0,\infty[\,\longrightarrow \mathbb{R}$ by
\begin{equation}\label{q0}
T_{\alpha}f(x):=
\lim_{\epsilon\rightarrow 0}
\frac{f(x+\epsilon x^{1-\alpha})-f(x)}{\epsilon},\quad x>0,
\end{equation}
provided the limit exists, in which case $f$ is called $\alpha$-differentiable. Additionally if $f$ is $\alpha$-differentiable on $\,]0,a[\,$ for some $a>0$, then $T_{\alpha}f(0):=\lim_{x\rightarrow 0^{+}}T_{\alpha}f(x)$. It is also shown
(\cite{khalil}, Theorem 2.2)
that if a function $f$ is differentiable at some $x>0$, then $T_{\alpha}f(x)$ exists and
\begin{equation}\label{q1}
T_{\alpha}f(x)=x^{1-\alpha}f^{\prime}(x).\\
\end{equation}
\indent It is pointed out in \cite{Ortigueira}
that the conformable derivative (\ref{q0})
lacks the properties (ii) and (iv).\\
\indent In this paper, we demonstrate that if $f$ is $\alpha$-differentiable
in the suggested conformable sense
at some $x>0$, then it must be differentiable (in the classical sense)
at $x$ and (\ref{q1}) is satisfied. This means
that differentiability and $\alpha$-differentiability
in the sense of the existence of the limit
(\ref{q0}) are in fact equivalent at any $x>0$. In particular,
a function $f:[0,\infty[\,\longrightarrow \mathbb{R}$ is $\alpha$-differentiable
at $x=0$ if and only if it is differentiable
on $\,]0,a[\,$ for some $a>0$, and
\begin{equation}\label{atz}
T_{\alpha}f(0)=\lim_{x\rightarrow 0^{+}}x^{1-\alpha}f^{\prime}(x).
\end{equation}
Thus, we have $T_{\alpha}f(0)=0$ for every function $f$
continuously differentiable on a neighbourhood of $x=0$
or, more generally, differentiable on a right neighbourhood of $x=0$ with a bounded derivative therein. \\
\indent The fact that the
conformable $\alpha$-derivative is not a fractional derivative was already implicitly pointed out by
Tarasov in \cite{Tarasov} where he proves
that the violation of the Leibniz rule,
$D^{\alpha}(fg)=f D^{\alpha}g+g D^{\alpha}f$,
is necessary for the order $\alpha$
of a differential operator $D^{\alpha}$
to be fractional. The conformable
operator $T_{\alpha}$ does satisfy the Leibniz rule
(\cite{khalil}, Theorem 2.2).\\
\indent First, restricting our attention to the conformable derivative helped us obtain a simpler proof than Tarasov's general result in \cite{Tarasov}. Second,
not only do we prove that the conformable derivative is not fractional, but we also express
$T_{\alpha}f$ explicitly in terms of $f^{\prime}$
by means of (\ref{q1}).
Using this explicit pointwise relation trivializes (and sometimes invalidates) many of the proofs/computational analysis in papers such as
\cite{Abdeljawad,Benaoumeur,Chen,Eslami,Ghazala,
Hosseini,Katugampola,khalil,Morales,
Ortigueira,Tarasov,Yang,Zhou} just to mention a few.
To see this, assume one adopts
the conformable derivative
 $T_{\alpha}$ to interpret fractional derivatives and considers, say, solving the fractional ODE $\,F(x,u,T^{\alpha}u,Du)=0$,
or the fractional PDE
$\,G(x,y,u,T^{\alpha}_{x}u,T^{\beta}_{y}u,D_{x}u,D_{y}u)=0$. Then, in reality, they are solving the ODE
$\,\widetilde{F}(x,u,Du)=0$, or the PDE
$\,\widetilde{G}(x,y,u,D_{x}u,D_{y}u)=0$, respectively, where
\begin{equation*}
\widetilde{F}(x,y,z):=F(x,y,x^{1-\alpha}z,z),
\;\widetilde{G}(x,y,z,w,v):=
G(x,y,z,x^{1-\alpha}w,y^{1-\beta}v,w,v).
\end{equation*}
\indent
We analyse a general fractional model of viscoelasticity and show that
the solutions that comes from interpreting
fractional derivatives in the conformable
sense are inconsistent with the ones
that correspond to the Riemann-Liouville and Caputo derivative. The latter are known to be in excellent fit with experimental studies. \\
\indent We conclude with the remodelling two real-world phenomena already modelled using the classical fractional derivative, the viscoelastic
deformation (the creep effect) in the tight sandstone 
H20-6 \cite{Ding}, and the world population growth
\cite{Almeida}, using the conformable derivative. 
We find that the error in the conformable derivative framework is larger than the error in the fractional framework by a considerable margin in both models.\\
\indent  The effect
of not using an adequate definition
for the fractional derivatives is not
usually addressed by researchers.
This study shows the importance of using a proper interpretation for fractional models.
\section{A misleading counterexample and incorrect statements}
\indent It is claimed in \cite{khalil} that
a function $\alpha$-differentiable in the conformable sense is not necessarily
differentiable. We shall prove in Theorem \ref{thm1}
below that this is incorrect at any $x>0$. The only counterexample
in the literature that allegedly backs up this claim is the example (see \cite{khalil}) of $g(x):=\sqrt{x}$.
Of course, while $T_{\frac{1}{2}}g(0)=1$,
$g^{\prime}(0)$ does not exist. First, this does not apply to the translation of $g$ to any $x_{0}>0$. Simply consider
$h(x):=\sqrt{x-x_{0}}$ with $x_{0}>0$.
Then $T_{\alpha}h(x_{0})={x_{0}^{\frac{1-\alpha}{2}}}
\lim_{\epsilon\rightarrow 0}\frac1{
\sqrt{\epsilon}}$
does not exist for any $0\leq \alpha<1$, and neither does $h^{\prime}(x_{0})$. Second,
since the existence of $f^{\prime}(0)$
is independent of the existence of
$\lim_{x\rightarrow 0^{+}} x^{1-\alpha}
f^{\prime}(x)$, we realize from (\ref{atz})
that the differentiability at $x=0$ is independent of
the conformable $\alpha$-differentiability there. For instance,
the function $\tilde{g}:=x^2
\chi_{\raisebox{-.5ex}{$\scriptstyle \mathbb{Q}$}}$, where $\chi_{\raisebox{-.5ex}{$\scriptstyle \mathbb{Q}$}}$ is the characteristic function of the rational numbers, is differentiable only at 0. Therefore,
by (\ref{atz}),
$T_{\alpha}\tilde{g}(0)$ does not exist.\\
\indent
Another analogous attempt to define a conformable
$\alpha$-derivative of a function
$f:[0,\infty[\,\rightarrow \mathbb{R}$
appears in
\cite{Katugampola1} where
(\ref{q0}) is reformulated to
\begin{equation}\label{ualpha}
U_{\alpha}f(x):=
\lim_{\epsilon\rightarrow 0}
\frac{f(x e^{\epsilon x^{-\alpha}})-f(x)}{\epsilon},\quad x>0,
\end{equation}
$U_{\alpha}f(0):=\lim_{x\rightarrow 0^{+}}
x^{1-\alpha}U_{\alpha}f(x)$, provided the respective limit exists.
It is shown in (\cite{Katugampola1}, Theorem 2.3)
that if $f^{\prime}(x)$ exists
at some $x>0$, then so does
$U_{\alpha}f(x)$ and
\begin{equation}\label{qq1}
U_{\alpha}f(x)=x^{1-\alpha}f^{\prime}(x).
\end{equation}
Again, it is claimed in \cite{Katugampola1}
that there exists a
conformable $\alpha$-differentiable
function that is not differentiable without providing
any proof. We shall show in Theorem \ref{thm2} in the next section that, with the exception of
the origin, this claim is also
incorrect.
\section{Main results}\label{sec3}
\begin{theorem}\label{thm1}
Fix $\,0<\alpha<1$ and let $x>0$. A function
$\,f:[0,\infty[\,\longrightarrow \mathbb{R}\,$
has a conformable fractional derivative of
order $\alpha$ at $x$ if and only if
it is differentiable at $x$ and
(\ref{q1}) holds.
\end{theorem}
\begin{proof}
By Theorem 2.2 in \cite{khalil}, it suffices
to prove that if $T_{\alpha}f(x)$
exists then so does $f^{\prime}(x)$.
We have
\begin{eqnarray*}
\lim_{\epsilon\rightarrow 0}
\frac{f(x+\epsilon)-f(x)}{\epsilon}
 &=&\lim_{\epsilon\rightarrow 0} x^{\alpha-1}
\frac{f(x+(\epsilon\,x^{\alpha-1})\,x^{1-\alpha})-f(x)}{
\epsilon\,x^{\alpha-1} }   \\
   &=&x^{\alpha-1}\lim_{\epsilon\rightarrow 0}\,
\frac{f(x+\epsilon\,x^{1-\alpha})-f(x)}{
\epsilon}.
\end{eqnarray*}
\end{proof}
\begin{theorem}\label{thm2}
Fix $\,0<\alpha<1$ and $x>0$. A function
$\,f:[0,\infty[\,\longrightarrow \mathbb{R}\,$
has a conformable fractional derivative $U_{\alpha}f(x)$ if and only if
it is differentiable at $x$ and
(\ref{qq1}) holds.
\end{theorem}
\begin{proof}
By Theorem 2.3 in \cite{Katugampola1}, we only need
to show that $f^{\prime}(x)$
exists whenever $U_{\alpha}f(x)$ does.
Since
\begin{equation*}
x+\epsilon=
e^{\log{\left(x+\epsilon\right)}}=
xe^{\log{\left(1+\frac{\epsilon}{x}\right)}}
= xe^{\epsilon
x^{-\alpha}
\left(x^{\alpha}\log{\left(1+\frac{\epsilon}{x}\right)^{
\frac{1}{\epsilon}}}\right)}
\end{equation*}
then
\begin{align*}
&\lim_{\epsilon\rightarrow 0}
\frac{f(x +\epsilon)-f(x)}{\epsilon}
 \;=\;\lim_{\epsilon\rightarrow 0}
\frac{f\left(xe^{\epsilon
x^{-\alpha}
\left(x^{\alpha}\log{\left(1+\frac{\epsilon}{x}\right)^{
\frac{1}{\epsilon}}}\right)}\right)-f(x)}{\epsilon}=\\
&=\;\lim_{\epsilon\rightarrow 0}
\frac{f\left(xe^{\epsilon
x^{-\alpha}
\left(x^{\alpha}\log{\left(1+\frac{\epsilon}{x}\right)^{
\frac{1}{\epsilon}}}\right)
}\right)-f(x)}{\epsilon\left(x^{\alpha}
\log{\left(1+\frac{\epsilon}{x}\right)^{
\frac{1}{\epsilon}}}\right)}\,\left(x^{
\alpha}\log{\left(1+\frac{\epsilon}{x}\right)^{
\frac{1}{\epsilon}}}\right).
\end{align*}
Observing that
\begin{equation*}
\lim_{\epsilon\rightarrow 0}
\log{\left(1+\frac{\epsilon}{x}\right)^{
\frac{1}{\epsilon}}}=\frac{1}{x}
\end{equation*}
we infer that
\begin{equation*}
\lim_{\epsilon\rightarrow 0}
\frac{f(x +\epsilon)-f(x)}{\epsilon}
\,=\,x^{\alpha-1}\lim_{\epsilon\rightarrow 0}
\frac{f\left(xe^{\epsilon
x^{-\alpha}}\right)-f(x)}{\epsilon}.
\end{equation*}
\end{proof}
We conclude with a general principle:
\begin{theorem}\label{genprin}
Suppose $h:\,]-1,1[\,\times\mathbb{R}\longrightarrow
\mathbb{R}$ is such that
$\,\lim_{\epsilon\rightarrow 0}h(\epsilon,x_{0})\neq 0$
for some $x_{0}\in \mathbb{R}$. Then a function $\phi:\mathbb{R}\longrightarrow
\mathbb{R}$ is differentiable at $x_{0}$
if and only if the limit
\begin{equation*}
\tilde{\phi}(x_{0}):=\lim_{\epsilon\rightarrow 0}
\frac{\phi\left(x_{0}+\epsilon h(\epsilon,x_{0})\right)-\phi(x_{0})}{\epsilon}
\end{equation*}
exists, in which case $\tilde{\phi}(x_{0})=
\rho_{1}(x_{0})\phi^{\prime}(x_{0})$ and
$\; \rho_{1}(x)=
\lim_{\epsilon\rightarrow 0}h(\epsilon,x)$.
\end{theorem}
\section{Applications}\label{secappl}
Theorems \ref{thm1} and \ref{genprin}
simply prove the pointlessness of the results in
\cite{khalil} and \cite{Abdeljawad}, respectively. Similarly, Theorem \ref{thm2} trivializes all
the results in \cite{Anderson,Katugampola1}.
For example, take Theorem 2.3 in
\cite{khalil}. It asserts
that if $f\in C([a,b])$ is
$\alpha$-differentiable on $]a,b[$, $a>0$,
and $f(a)=f(b)$, then there exists $c\in\,]a,b[\,$
such that $T_{\alpha}f(c)=0$. But this is immediate
from the classical Rolle's theorem. Indeed,
by Theorem \ref{thm1}, the function $f$ is
differentiable on $]a,b[$. So, there exists
$c\in\,]a,b[\,$ such that $f^{\prime}(c)=0$. Using Theorem 1 again implies $T_{\alpha}f(c)=c^{1-\alpha}f^{\prime}(c)=0$.
Next, Theorem 2.4 in \cite{khalil} argues
that if $f\in C[a,b]$ is $\alpha$-differentiable
on $]a,b[$, then there is $c\in\,]a,b[$ such that
$T_{\alpha}f(c)=\frac{f(b)-f(a)}{
\frac{b^{\alpha}}{\alpha}-\frac{a^{\alpha}}{\alpha}}$.
Again, by Theorem \ref{thm1}, we know
that $f$ is differentiable on $]a,b[$ and
$T_{\alpha}f(c)=\frac{f^{\prime}(c)}{c^{\alpha-1}}$.
Therefore, the claim results from applying the classical mean value property to the functions $f$ and
$x\mapsto x^{\alpha}/\alpha$, and so on for the rest of \cite{khalil}.\\
\indent
For another example, consider the results
in \cite{Abdeljawad}. They are based on Definition 2.1 in the same reference that introduces a right (left) conformable $\alpha$-derivative, $T^{a}_{\alpha}f$ ($^{b}_{\alpha}Tf$),
where $f$ is defined on the right of $x=a$
(on the left of $x=b$) through the limits
\begin{align*}
   T^{a}_{\alpha}f(x)&=\lim_{\epsilon\rightarrow 0}
\frac{f(x+\epsilon (x-a)^{1-\alpha})-f(x)}{\epsilon},\quad x>a, \\
   ^{b}_{\alpha}Tf(x)&=
\lim_{\epsilon\rightarrow 0}
\frac{f(x+\epsilon (b-x)^{1-\alpha})-f(x)}{\epsilon},\quad x<b.
\end{align*}
Consider a continuous function $f:\,]a,b[\,\longrightarrow \mathbb{R}$
and let $a<x<b$. Applying Theorem \ref{genprin}, we discover that
the following are equivalent:
\begin{itemize}
  \item $f$ is right $\alpha$-differentiable at $x$
  \item $f$ is left $\alpha$-differentiable at $x$
  \item $f$ is differentiable at $x$
\end{itemize}
Furthermore, we find from Theorem \ref{genprin} that
we do not need to assume the differentiability of $f$
to apply the relations
\begin{equation*}
T^{a}_{\alpha}f(x)=
(x-a)^{1-\alpha}f^{\prime}(x),\; x>a,\quad
^{b}_{\alpha}Tf(x)=
(b-x)^{1-\alpha}f^{\prime}(x),\; x<b.
\end{equation*}
These follow directly from the $\alpha$-differentiability assumption that guarantees
the differentiability.
Taking this into account, one sees how the results in \cite{Abdeljawad} follow directly from classical calculus. \\
\indent All the results in \cite{Anderson} are built on the premise of an inaccurate citation (see \cite{Anderson}, part II) of item (7) of
Theorem 2.3 in \cite{Katugampola1} which reads
\begin{theorem}(\cite{Katugampola1}, Theorem 2.3 (7))
Let $\alpha\in\, ]0,1]$ and $f$ be $\alpha$-differentiable at a point $x>0$;
so that the limit (\ref{ualpha}) exists.
\textbf{If}, \textbf{in addition},
$\textbf{f}$ \textbf{is differentiable},
\textbf{then} (\ref{qq1}) holds.
\end{theorem}
The relation (\ref{qq1}) is used throughout \cite{Anderson} without imposing
the differentiability assumption. If differentiability
is assumed implicitly, then what is the novelty in studying the operator $x^{1-\alpha} D$ on differentiable functions? We have demonstrated in
Theorem \ref{thm2} that there is no
function $\alpha$-differentiable in the sense of
Katugampola \cite{Katugampola1} that is not differentiable.\\
\indent Interpreting fractional derivatives
in the conformable sense
alters fractional differential problems
into differential problems with the usual classical derivatives. A persistent question here is whether
the resulting differential equation still
accurately models the original fractional physical phenomena. Recently, some fractional models
were analysed after rewriting fractional derivatives in terms of the classic derivatives
via (\ref{q1}) without addressing this crucial
issue (see e.g., \cite{Chen,Eslami,Ghazala,
Hosseini,Morales,Yang,Zhou}). In fact, in many cases, these new equations happen to be already intensively studied.\\
\textbf{Example 1.}\\
\indent In \cite{Benaoumeur}, the authors seek the existence of a solution to the local fractional problem
\begin{equation*}
\left\{
  \begin{array}{ll}
   x^{(\alpha)}(t)=f(t,x(t)), & \hbox{$t\in [a,b],\;a >0$,} \\
 x(a)=x_{0},&
  \end{array}
\right.
\end{equation*}
where $f:[a,b]\times \mathbb{R}\rightarrow \mathbb{R}$
is a continuous function. Since the interpretation (\ref{q0}) for fractional derivatives is adopted in \cite{Benaoumeur}, then
what is really studied therein is the existence for the very simple well-understood local IVP
\begin{equation*}
\left\{
  \begin{array}{ll}
   x^{\prime}(t)=g(t,x(t)), & \hbox{$t\in [a,b],\;a >0$,} \\
 x(a)=x_{0},&
  \end{array}
\right.
\end{equation*}
where $g(t,s):=s^{\alpha-1}f(t,s)$. Observe that $g$ is also continuous, and if $f$ is Lipschitz, then so is $g$.\\
\textbf{Example 2.}\\
\indent Applying Theorem \ref{thm1}, one immediately
realizes that the fractional differential equations
\begin{align*}
  &y^{(\frac{1}{2})}+y=x^2+2x^{\frac{3}{2}},\; y(0)=0,\\
   & y^{(\alpha)}+y=0,\; 0<\alpha\leq 1\\
&y^{(\frac{1}{2})}+\sqrt{x} y=xe^{x}, \\
&y^{(\frac{1}{2})} = \frac{x^{\frac{3}{2}}+y\sqrt{x}}{2x+3y}
\end{align*}
are not really solved in \cite{khalil}.
What are really solved there are
the respective elementary ODEs:
\begin{align*}
  & y^{\prime}+x^{-\frac{1}{2}}y=x^{\frac{3}{2}}+2x,\; y(0)=0,\\
   & x^{1-\alpha}y^{\prime}+y=0,\\
& y^{\prime}+ y=\sqrt{x}e^{x}, \\
& y^{\prime} = \frac{x+y}{2x+3y}.
\end{align*}
\section{Comparative analysis of a fractional model}\label{companal}
\begin{center}
\begin{figure}[h]
  \centering
  \includegraphics[width=9 cm]{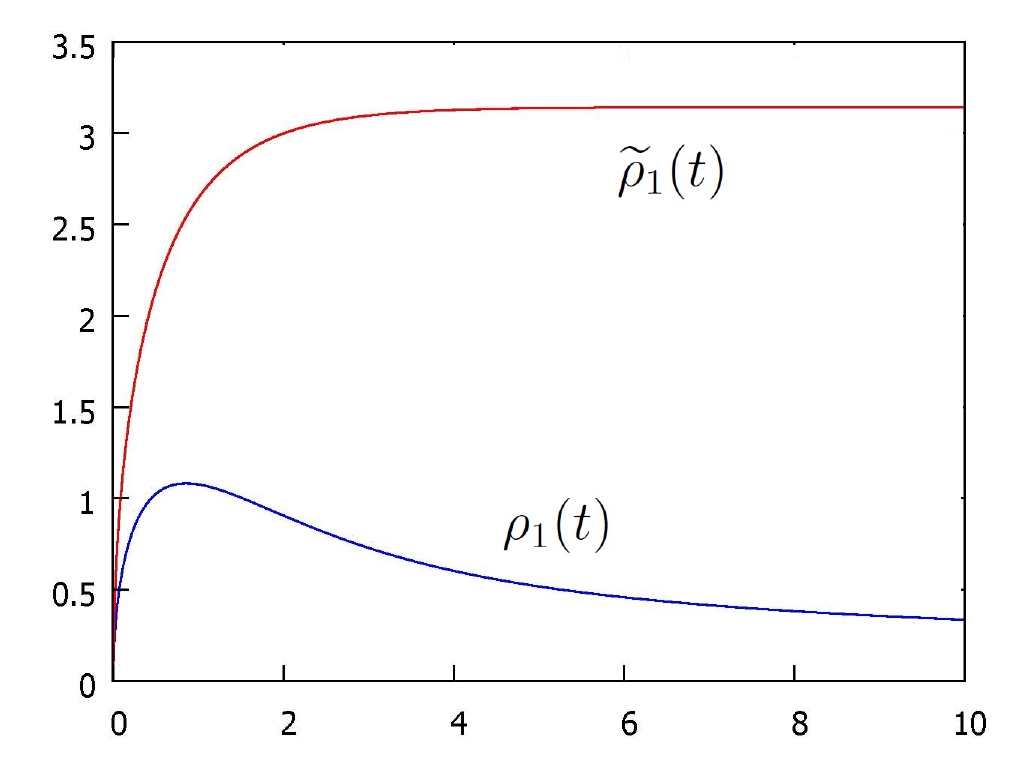}\\
  \caption{The profile of the strain
that corresponds to an exponentially decaying stress
on a viscoelastic material}
\label{fig1}
\end{figure}
\end{center}

\begin{center}
\begin{figure}[h]
  \centering
  \includegraphics[width=9 cm]{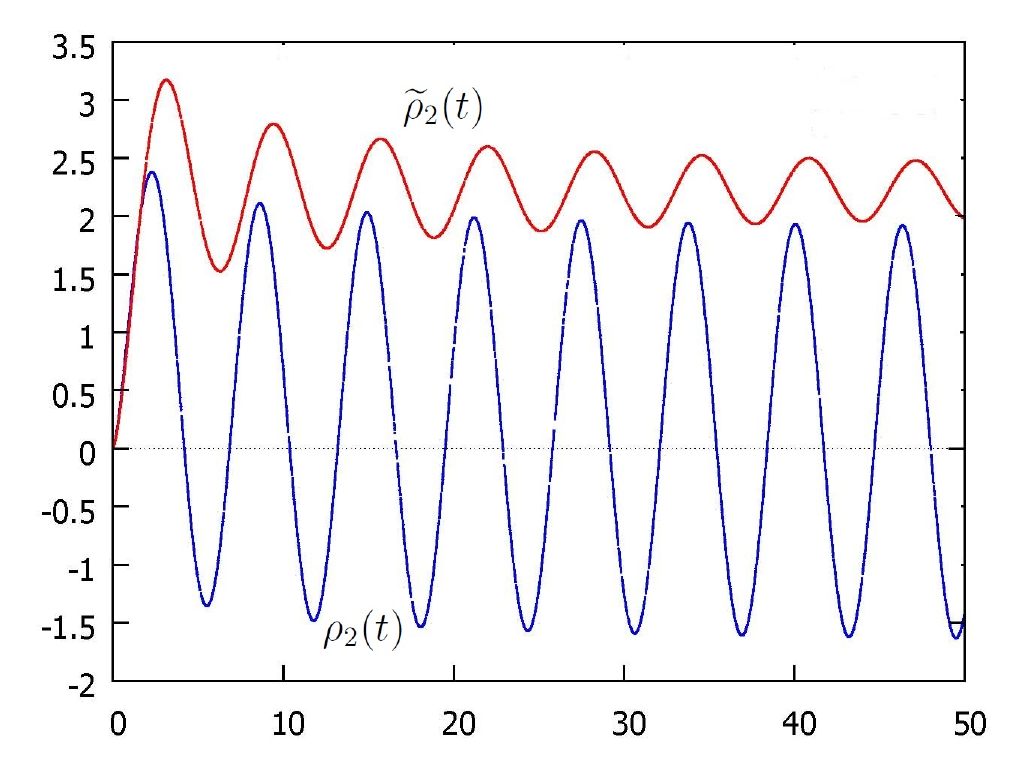}\\
  \caption{The profile of the strain
that corresponds to an oscillating stress
on a viscoelastic material ($A_{1}=1$).}
\label{fig2}
\end{figure}
\end{center}

\begin{center}
\begin{figure}[h]
  \centering
  \includegraphics[width=9 cm]{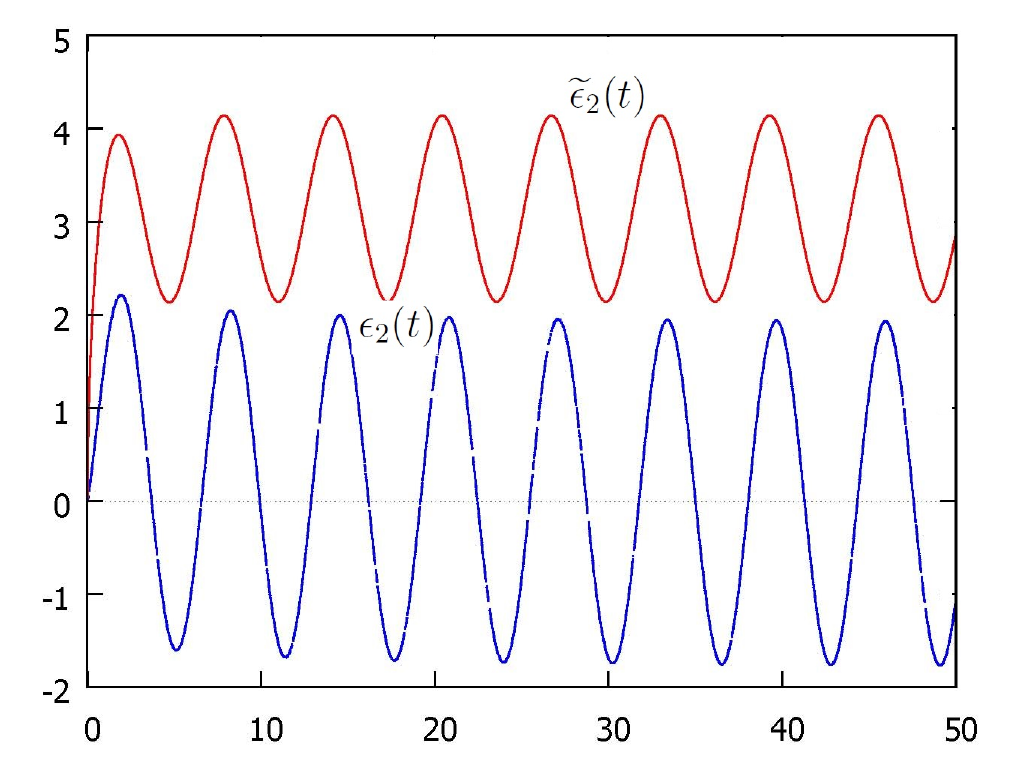}\\
  \caption{The strain
that corresponds to an oscillating stress
on a viscoelastic material ($b=A_{1}=1$).}
\label{fig3}
\end{figure}
\end{center}

Caputo and Mainardi proposed in \cite{caputom}
the following fractional model for the viscoelastic
reaction:
\begin{equation}\label{model0}
\sigma+b D^{\alpha}\sigma=
A_{0}\epsilon+A_{1}D^{\alpha}\epsilon,
\end{equation}
where $\sigma$ is the stress (tension),
$\epsilon$ is the strain (deformation), and
$b$, $A_{0}$, $A_{1}$ are constants.
The model (\ref{model0}) was found to be in
agreement with experimental
data of various materials. Later,
Bagley and Torvik investigated this model
in a series of studies that promoted it as
the definitive model. They showed in \cite{bagley} that, with the restrictions
\begin{equation*}
b\geq 0,\; A_{0}\geq 0,\; A_{1}>0,\; A_{1}\geq A_{0}b,
\end{equation*}
and with $D^{\alpha}$ taken to be the Riemann-Liouville
or Caputo fractional derivative, the model (\ref{model0})
fulfills the laws of thermodynamics.
When $D^{\alpha}$ is interpreted as
the Riemann-Liouville derivative, the general solution of
(\ref{model0}) is given by:
\begin{equation}\label{RLsol}
\epsilon(t)=
\frac{b}{A_{1}}\sigma(t)
-C_{1}\frac{
E_{\alpha,\alpha}\left(-\frac{A_{0}}{A_{1}}(t-a)
^{\alpha}\right)}{(t-a)^{1-\alpha}}+\frac{1-bA_{0}}{A_{1}^{2}}
\int_{a}^{t}
\frac{E_{\alpha,\alpha}\left(-\frac{A_{0}}{A_{1}}(t-\xi)
^{\alpha}\right)}{(t-\xi)^{1-\alpha}}\sigma{(\xi)}d\xi
\end{equation}
where
\begin{equation*}
C_{1}=\Gamma(\alpha)\lim_{t\rightarrow a^{+}}
\left(t-a\right)^{1-\alpha}
\left( \frac{b}{A_{1}}\sigma(t)-\epsilon(t)\right)
\end{equation*}
and
\begin{equation*}
E_{\alpha,\beta}(z)=\sum_{k=0}^{\infty}\frac{z^k}{
\Gamma\left(\alpha k+\beta\right)},\;\;
\beta, z\in \mathbb{C}, \alpha>0
\end{equation*}
is the generalized Mittag-Leffler function \cite{Kilbas}.
If $D^{\alpha}$ is the Caputo derivative, then
the general solution of
(\ref{model0}) is:
\begin{equation}\label{Csol}
\epsilon(t)=
\frac{b}{A_{1}}\sigma(t)
+
\int_{a}^{t}
\frac{E_{\alpha,\alpha}\left(-\frac{A_{0}}{A_{1}}(t-\xi)
^{\alpha}\right)}{(t-\xi)^{1-\alpha}}
\left( \frac{1-bA_{0}}{A_{1}^{2}}\sigma(\xi)+C_{2}\frac
{A_{0}}{A_{1}}\right)
d\xi+C_{2}
\end{equation}
where
\begin{equation*}
C_{2}=\left( \frac{b}{A_{1}}\sigma(a)-\epsilon(a)\right).
\end{equation*}
Notice that, if $\sigma$ and $\epsilon$ are right continuous at $t=a$ and
\begin{equation}\label{initial}
b\sigma(a)-A_{1}\epsilon(a)=0,
\end{equation}
then $C_{1}=C_{2}=0$, and the solutions
(\ref{RLsol}) and (\ref{Csol}) coincide and take the form
\begin{equation}\label{sol}
\epsilon(t)=
\frac{b}{A_{1}}\sigma(t)
+\frac{1-bA_{0}}{A_{1}^{2}}
\int_{a}^{t}
E_{\alpha,\alpha}\left(-\frac{A_{0}}{A_{1}}(t-\xi)
^{\alpha}\right)
\frac{\sigma(\xi)}{(t-\xi)^{1-\alpha}}d\xi.
\end{equation}
We hereafter consider this case.\\
\indent
Now, if we use the conformable derivative
(\ref{q0}), then, using Theorem \ref{thm1},
the model (\ref{model0}) would attain the general solution
\begin{equation*}
\widetilde{\epsilon}(t)=
\frac{b}{A_{1}}\sigma(t)-\frac{b}{A_{1}}\sigma(a)
e^{\frac{A_{0}}{A_{1}}\frac{a^{\alpha}-t^{\alpha}}{\alpha}}+
C_{3}e^{-\frac{A_{0}}{A_{1}}\frac{t^{\alpha}}{\alpha}}+
\frac{A_{1}-bA_{0}}{A^{2}_{1}}
e^{-\frac{A_{0}}{A_{1}}\frac{t^{\alpha}}{\alpha}}
\int_{a}^{t}\frac{\sigma(\xi)}{\xi^{1-\alpha}}
e^{\frac{A_{0}}{A_{1}}\frac{\xi^{\alpha}}{\alpha}}d\xi,
\end{equation*}
which solves the ODE:
\begin{equation*}
A_{1}t^{1-\alpha}\widetilde{\epsilon}^{\prime}+A_{0}\widetilde{\epsilon}=
b t^{1-\alpha}\sigma^{\prime}+\sigma,\;
t\geq a.
\end{equation*}
If $\sigma e^{\frac{A_{0}}{A_{1}}\frac{\xi^{\alpha}}{\alpha}}/\xi^{1-\alpha}$ is integrable and the initial condition (\ref{initial}) is satisfied, then
\begin{equation}\label{conformsol}
\widetilde{\epsilon}(t)=
\frac{b}{A_{1}}\sigma(t)+
\frac{A_{1}-bA_{0}}{A^{2}_{1}}
e^{-\frac{A_{0}}{A_{1}}\frac{t^{\alpha}}{\alpha}}
\int_{a}^{t}\frac{\sigma(\xi)}{\xi^{1-\alpha}}
e^{\frac{A_{0}}{A_{1}}\frac{\xi^{\alpha}}{\alpha}}d\xi.
\end{equation}
Let us compare $\widetilde{\epsilon}$
to the solution $\epsilon$ in (\ref{sol})
for two cases:
\begin{enumerate}
  \item the stress decays exponentially:
$\;\sigma_{1}(t)=e^{-t}$
  \item the stress is a wave:
$\;\sigma_{2}(t)=\sin{t}$
\end{enumerate}
A detailed
practical background on the decaying
and oscillatory stress
can be found in \cite{Lakes} and \cite{Ward}.\\
\indent Let $\epsilon_{i}, \widetilde{\epsilon}_{i}$ denote
the strain that corresponds
to the stress $\sigma_{i}$, $i=1,2$, and
let $a=0$, $\alpha=0.5$. For simplicity, we study the model (\ref{model0}) with $A_{0}=0$.\\
\indent
The term $\frac{b}{A_{1}}\sigma$, being the first one in both expressions (\ref{sol}) and (\ref{conformsol}),
is obviously irrelevant in our comparison. So, we
compare instead $\rho_{i}$ to $\widetilde{\rho}_{i}$
where
\begin{equation*}
\rho_{i}(t) :=A^{2}_{1}\left(\epsilon_{i}(t)-\frac{b}{A_{1}}
\sigma_{i}(t)\right),\quad
\widetilde{\rho}_{i}(t):=A^{2}_{1}\left(\widetilde{\epsilon}_{i}(t)-
\frac{b}{A_{1}}\sigma_{i}(t)\right).
 \end{equation*}
The functions $\rho_{i}$ and $\widetilde{\rho}_{i}$
have the same features as
$\epsilon_{i}$ and $\widetilde{\epsilon}_{i}$, respectively.\\
\indent Consider the stress $\sigma_{1}$. Then
\begin{align*}
\rho_{1}(t)&=\,\frac{1}{\sqrt{\pi}}\int_{0}^{t}\frac{e^{-\xi}}{\sqrt{t-\xi}}
d\xi=
\mathrm{\emph{erfi}}{\left(\sqrt{t}\right)
}e^{-t},\\
\widetilde{\rho}_{1}(t)&=\,\frac{A_{1}}{\sqrt{\pi}}
\int_{0}^{t}\frac{e^{-\xi}}{\sqrt{\xi}}d\xi=A_{1} \mathrm{\emph{erf}}{\left(\sqrt{t}\right)},
\end{align*}
where $\mathrm{\emph{erf}}$ and
$\mathrm{\emph{erfi}}$ are the error and imaginary error functions, respectively.\\
\\ By the properties of the error function, $\widetilde{\rho}_{1}$
is strictly increasing with $\lim_{t\longrightarrow\infty}\widetilde{\rho}_{1}(t)=
 A_{1}$. This implies that
\begin{equation*}
\lim_{t\longrightarrow\infty}\widetilde{\epsilon}_{1}(t)=
\frac{\pi}{A_{1}}.
\end{equation*}
Thus, the conformable interpretation
does not show the elasticity effect, as it implies
a permanent deformation of the material despite
the rapid decrease of the stress.
\\ \indent To understand the behavior of $\rho_{1}$,
write $\psi(t)=\sqrt{\pi}\rho_{1}(t)$. Then
\begin{equation*}
\psi(t)=
\int_{0}^{t}\frac{e^{-\xi}}{\sqrt{t-\xi}}d\xi=
e^{-t}\int_{0}^{t}\frac{e^{\xi}}{\sqrt{\xi}}d\xi.
\end{equation*}
Since $\frac{e^{\xi}}{\sqrt{\xi}}\in L^{1}([0,t])$
for all $t<\infty$, then we may calculate
\begin{equation*}
\lim_{t\rightarrow \infty}\psi(t)=
\lim_{t\rightarrow \infty}\frac{
\int_{0}^{t}\frac{e^{\xi}}{\sqrt{\xi}}d\xi}{e^{t}}=
\lim_{t\rightarrow \infty}\frac{1}{\sqrt{t}}=0.
\end{equation*}
Moreover,
\begin{equation*}
\psi^{\prime}(t)=\frac{\phi(t)}{e^{t}},\;
\phi(t):=
\frac{e^{t}}{\sqrt{t}}-\int_{0}^{t}\frac{
e^{\xi}}{\sqrt{\xi}}d\xi.
\end{equation*}
Since $\phi(1)=e-\sqrt{\pi}\mathrm{\emph{erfi}}(1)<0$
and $\phi^{\prime}(t):=
-\frac{e^{t}}{2\sqrt{t^3}}<0$,
then $\phi<0$
on $[1,\infty[$. This shows that $\rho_{1}^{\prime}<0$
on $[1,\infty[$. Consequently, ${\rho}_{1}$ is strictly decreasing on $[1,\infty[$, and we have $\lim_{t\longrightarrow\infty}{\rho}_{1}(t)=0$.
This is in line with the elasticity effect
as the strain decreases with the decrease in the stress (see Figure \ref{fig1}).
\indent Now, consider the stress $\sigma_{2}$. We have
\begin{align*}
\rho_{2}(t)&=\,\frac{1}{\sqrt{\pi}}
\int_{0}^{t}\frac{\sin{\xi}}{\sqrt{t-\xi}}d\xi=
\sqrt{2}\left(
\mathrm{\emph{C}}\left( \sqrt{\frac{2}{\pi}}\sqrt{t}\right)\sin{t}-
\mathrm{\emph{S}}\left( \sqrt{\frac{2}{\pi}}\sqrt{t}\right)\cos{t}\right),\\
\widetilde{\rho}_{2}(t)&=\,\frac{A_{1}}{\sqrt{\pi}}
\int_{0}^{t}\frac{\sin{\xi}}{\sqrt{\xi}}d\xi=
\sqrt{2}A_{1}\mathrm{\emph{S}}\left( \sqrt{\frac{2}{\pi}}\sqrt{t}\right),
\end{align*}
where
$\mathrm{\emph{S}}$ and $\mathrm{\emph{C}}$ are, respectively, the sine and cosine Fresnel integrals (see Figure \ref{fig2}).
Let us further investigate
\begin{equation*}
\widetilde{\epsilon}_{2}(t)=
\frac{1}{A_{1}}
\left(\frac{\widetilde{\rho}_{2}(t)}{A_{1}}
+b\sigma_{2}(t) \right)=
\frac{1}{A_{1}}
\left(\int_{0}^{t}\frac{\sin{\xi}}{\sqrt{\xi}}d\xi
+b\sigma_{2}(t) \right).
\end{equation*}
\indent There exists $T>0$ such that
$\displaystyle \int_{0}^{t}\frac{\sin{\xi}}{\sqrt{\xi}}d\xi
+\sin{t}>0$ for all $t\geq T$. This is because
\begin{equation*}
\lim_{t\rightarrow\infty}
\int_{0}^{t}\frac{\sin{\xi}}{\sqrt{\xi}}d\xi=
\int_{0}^{\infty}\frac{\sin{\xi}}{\sqrt{\xi}}d\xi=
\sqrt{\frac{\pi}{2}}>1.
\end{equation*}
So, for every $\delta>0$, there exists $T>0$ such that
\begin{equation*}
\left|\int_{0}^{t}\frac{\sin{\xi}}{\sqrt{\xi}}d\xi
-\sqrt{\frac{\pi}{2}}\right|<\delta
\end{equation*}
for all $t\geq T$. \\
\indent Therefore, if the material is such that $b\leq 1$, or
the input stress $\sigma_{2}$ is updated to be
$\sin{t}/b$, then $\widetilde{\epsilon}_{2}(t)>0$.
This means a permanent deformation after some time. The strain $\epsilon_{2}$, on the other hand, may oscillate, i.e., change sign. Once again, we see how
the conformable derivative
leads to a misunderstanding of the elasticity
nature in certain viscoelastic materials
(see Figure \ref{fig3}).
\section{Efficiency in modelling a real-world problem}
\subsection{Gross National Production}
*********
R. Almeida, What is the best fractional derivative to fit data?. Appl. Anal. Discrete
Math. 11, No 2 (2017), 358-368.

\subsection{World population}
A fractional model for the world population growth between the years 1910 and 2010 is proposed in \cite{Almeida}:
\begin{equation}\label{fracpop}
D^{\alpha}N(t)=PN_{t}.
\end{equation}
Here, $N(t)$ is the population number at the time
$t$ and $P$ is the growth rate. The fractional derivative $D^{\alpha}$ is taken in the Caputo sense. The solution of (\ref{fracpop}) then takes the form
\begin{equation}\label{solfraccap}
N(t)=N_{0}E_{\alpha,1}(P t^{\alpha}),
\end{equation}
where $N_{0}$ is the population number at the initial time. The authors in \cite{Almeida}
regard the fractional order $\alpha$
and the growth rate $P$ in (\ref{solfraccap}) as the parameters of the problem. They apply the sum of squares minimization
technique to determine these parameters such that
the function (\ref{solfraccap}) best fits
the census data provided by the united nations
\cite{un}.
The best parameters values turn out to be
$\alpha\approx 1.393298755$ and $P\approx 0.0034399$.
\\
\indent This method proved successful as the function
(\ref{solfraccap}) gives the population number
with an error of $E_{fractional}=2.050575\times 10^{5}$
as opposed to the classical integer-order model
$N^{\prime}(t)=PN(t)$ that yields the error
$E_{classical}=7.0795\times 10^{5}$. The errors represent the sum of the squares
of the differences between the data and the corresponding values of $N(t)$ in each case.\\
\indent If we reconsider the fractional model (\ref{fracpop}) and interpret $D^{\alpha}$
in the conformable sense, we get the population growth function
\begin{equation}\label{solpopcon}
\widetilde{N}(t)=N_{0}e^{p{t^{\alpha}}/{\alpha}}.
\end{equation}
To make a fair comparison, we follow the method proposed in \cite{Almeida} to determine the parameters $\alpha$
and $P$. First, we linearize (\ref{solpopcon}) to take the form $y=Ax+b$ where
$x=\log{t}$, $y=\log{\log{\left({\widetilde{N}(t)}/{N_{0}}\right)}}$,
$A=\alpha$, and $B=\log{\left({p}/{\alpha}\right)}$.
The least squares error minimization gives
$\alpha= 1.368138029$ and $P=0.003599846$.
Calculating the error the same way as $E_{fractional}$, we discover that $E_{conformable}=439528\times 10^{5}$ which is better than $E_{classical}$, yet
worse than $E_{fractional}$ by approximately the same margin. The figure below
is an update of Figure 2 in \cite{Almeida} after introducing the function $\widetilde{N}$.
\begin{center}
\begin{figure}[h]
  \centering
  \includegraphics[width=12 cm]{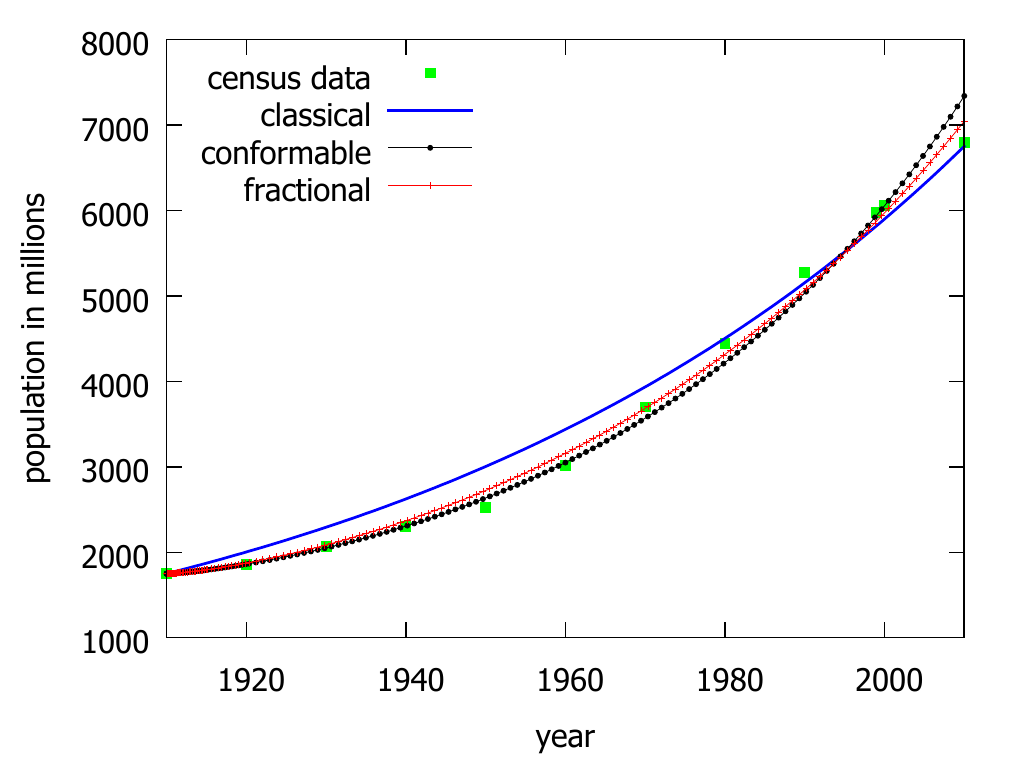}\\
\end{figure}
\end{center}
\section*{Acknowledgement}
The authors are grateful for
the valuable comments of the anonymous referees
that improved the presentation of this manuscript.
\section{Conclusion}
A definition for the fractional differential operator that allows the equivalence between the existence of the fractional derivative and differentiability might be acceptable. But, this is not the only issue with the conformable derivative. The more serious problem is that the conformable derivative is precisely the classical derivative times a root function. As tempting as it may be to express the fractional derivative in terms of the integer-order derivative, as the conformable interpretation suggests, it can break the fractional model and produce an irrelevant integer-order problem.
\section{Compliance with Ethical Standards}
$\quad$
The authors confirm that this work is original and has not been published elsewhere, nor is it currently under consideration for publication elsewhere.\\
\indent
We also declare that there is no conflict of interest, and that the study does not involve any human participants or animals and does not require
any form of a consent.

\end{document}